\documentclass{amsart}
\usepackage{amssymb}
\newtheorem{lemma}{Lemma}[section]
\newtheorem{propos}[lemma]{Proposition}
\newtheorem{theorem}[lemma]{Theorem}
\newtheorem{corol}[lemma]{Corollary}
\newtheorem{example}[lemma]{Example}

\newcommand{\CQ}{\hbox{{$\mathcal Q$}}}
\newcommand{\CS}{\hbox{{$\mathcal S$}}}
\newcommand{\CE}{\hbox{{$\mathcal E$}}}
\newcommand{\CC}{\hbox{{$\mathcal C$}}}
\newcommand{\CA}{\hbox{{$\mathcal A$}}}
\newcommand{\CH}{\hbox{{$\mathcal H$}}}
\newcommand{\CR}{\hbox{{$\mathcal R$}}}
\newcommand{\CL}{\hbox{{$\mathcal L$}}}
\newcommand{\cg}{\mathfrak{g}}

\newcommand{\C}{{\mathbb{C}}}
\newcommand{\R}{{\mathbb{R}}}

\newcommand{\Z}{{\mathbb{Z}}}

\newcommand{\del}{\partial}

 \renewcommand{\span}{{\rm span}}

\newcommand{\extd}{{\rm d}}
\newcommand{\trace}{{\rm tr}}
\newcommand{\isom}{{\cong}}
\newcommand{\eps}{{\epsilon}}
\newcommand{\tens}{\mathop{\otimes}}
\newcommand{\la}{{\triangleright}}

\newcommand{\Ad}{{\rm Ad}}
\newcommand{\id}{{\rm id}}
\newcommand{\<}{\langle}
\renewcommand{\>}{\rangle}
\newcommand{\End}{{\rm End}}

\newcommand{\Dsl}{D\!\!\!\!/}

\renewcommand{\o}{{}_{\scriptscriptstyle(1)}}
\renewcommand{\t}{{}_{\scriptscriptstyle(2)}}
\newcommand{\thr}{{}_{\scriptscriptstyle(3)}}

\newcommand{\bo}{{}^{\scriptscriptstyle\bar{(1)}}}
\newcommand{\bt}{{}^{\scriptscriptstyle\bar{(\infty)}}}
\newcommand{\bz}{{}^{\scriptscriptstyle\bar{(0)}}}

\newcommand{\rcocross}{{\blacktriangleright\kern -2pt<}}
\newcommand{\lcross}{{>\kern-4pt\triangleleft}}
\newcommand{\lcocross}{{>\kern-2pt\blacktriangleleft}}
\newcommand{\cobicross}{{\triangleright\kern-1pt\blacktriangleleft}}
\newcommand{\codcross}{{\blacktriangleright\kern-2pt\blacktriangleleft}}

\newcommand{\eproof}{$\quad \diamond$\goodbreak}

\newcommand{\eqn}[2]{\begin{equation}#2\label{#1}\end{equation}}

\begin{document}

\title[DIFFERENTIALS ON QUANTUM DOUBLES]{\rm \large CLASSIFICATION OF DIFFERENTIALS ON QUANTUM\\
DOUBLES AND FINITE NONCOMMUTATIVE GEOMETRY}
\author{S. Majid}
\address{School of Mathematical Sciences\\
Queen Mary, University of London\\ 327 Mile End Rd,  London E1
4NS, UK.}

\thanks{This paper is in final form and no version of it will be submitted for publication
elsewhere} 

\date{5/2002: resubmitted 1/2003}%

\maketitle

\begin{abstract} We discuss the construction of finite
noncommutative geometries on Hopf algebras and finite groups in
the `quantum groups approach'. We apply the author's previous
classification theorem, implying that calculi in the factorisable
case correspond to blocks in the dual, to classify differential
calculi on the quantum codouble $D^*(G)=kG\lcocross k(G)$ of a
finite group $G$. We give $D^*(S_3)$ as an example including its
exterior algebra and lower cohomology. We also study the calculus
on $D^*(\CA)$ induced from one on a general Hopf algebra $\CA$ in
general and specialise to $D^*(G)=U(\cg)\lcocross k[G]$ as a
noncommutative isometry group of an enveloping algebra $U(\cg)$ as
a noncommutative space.
\end{abstract}

\section{Introduction}

Coming out of the theory of quantum groups has emerged an approach
to noncommutative geometry somewhat different from the operator
algebras and K-theory one of Connes and others but with some
points of contact with that as well as considerable contact with
abstract Hopf algebra theory. The approach is based on building up
the different layers of geometry: the differential structure, line
bundles, frame bundles, etc. eventually arriving at spinors and a
Dirac operator as naturally constructed and not axiomatically
imposed as a definition of the geometry (as in Connes' spectral
triple theory).

This article has three goals. The first, covered in Section~2, is
an exposition of the overall dictionary as well as the differences
between the Hopf algebraic quantum groups approach and the more well-known
operator algebras and K-theory one. We also discuss issues that appear to have
led to confusion in the literature, notably the role of the Dirac
operator. The quantum groups methods are very algebraically
computable and hence should be interesting even to readers coming
from the Connes spectral triple side of noncommutative geometry. In fact,
the convergence and interaction between the two approaches is a
very important recent development and we aim to bridge between
them with our overview.

The second aim, in Sections~3, 4, is a self-contained
demonstration of the starting point of the quantum groups
approach, namely the classification of bicovariant differential
calculi on Hopf algebras. We start with the seminal work of
Woronowicz \cite{Wor:dif} for any Hopf algebra, which is mostly
known to Hopf algebraists and is included mainly for completeness.
We then explain our classification\cite{Ma:cla} for factorisable
quantum groups and apply it to finite-dimensional Hopf algebras
such as $\C_q[SL_2]$ at a cube root of unity. We also discuss
several subtleties not so well known in the literature. Let us
mention that in addition there are the Beggs--Majid classification
theorem \cite{BegMa:dif} for bicrossproduct quantum groups, see
also \cite{NME:cla}, and the Majid--Oeckl twisting
theorem\cite{MaOec:twi} (which covers most triangular quantum
groups), which altogether pretty much cover all classes of
interest. In this way the first layer of geometry, choosing the
differential structure, is more or less well understood in the
quantum groups approach. The theory is interesting even for
finite-dimensional Hopf algebras as `finite geometries' as our
example of $\C_q[SL_2]$ demonstrates.

Our third aim is to present some specific new results about
differential calculi on quantum doubles. Section~5 applies the
classification theorem of Section~4 to differential calculi on
$D^*(\CA)$ where $\CA$ is a finite-dimensional Hopf algebra and
$D^*(\CA)$ is the dual of its Drinfeld double\cite{Dri}. We
demonstrate this on the nice example of the dual of the Drinfeld
quantum double of a finite group, which is in fact a cross
coproduct $kG\lcocross k(G)$. Even this case is interesting as we
demonstrate for $G=S_3$, where we compute also its exterior
algebra and cohomology. Section~6 looks at the general case of
$D^*(\CA)$, where we show that a calculus on $\CA$ extends in a
natural way to one on the double, generalising the finite group
case in \cite{NME:cla}. The results also make sense for
infinite-dimensional quantum doubles such as $D^*(G)$ in the form
$U(g)\lcocross k[G]$ for an algebraic group of Lie type. We
propose the general construction as a step towards a duality for
quantum differentials inspired by T-duality in physics.

\section{Comparing the K-theoretic and quantum
groups approaches}

A dictionary comparing the two approaches appears in Table~1, and
let me say right away that we are not comparing like with like.
The Connes K-theoretic approach is mathematically purer and has deeper
theorems, while the quantum groups one is more organic and
experience-led. Both approaches of course start with an algebra
$M$, say with unit (to keep things simple), thought of as
`functions' but not necessarily commutative. The algebra plays the
role of a topological space in view of the Gelfand-Naimark
theorem. I suppose this idea goes back many decades to the birth
of quantum mechanics. Integration of course is some kind of linear
functional. In either approach it has to be made precise using
analysis (and we do not discuss this here) but at a conceptual
level there is also the question: which functional? The cyclic
cohomology approach gives good motivation to take here an element
of $HC^0(M)$,
i.e. a trace functional as an axiom. We don't have a good axiom in
the quantum groups approach except when $M$ is itself a quantum
group, when translation invariance implies a unique one. In
examples, it isn't usually a trace. For other spaces we would hope
that the geometry of the situation, such as a (quantum group)
symmetry would give the natural choice. So here we already see a
difference in scope and style of the approaches. Also, since we
usually proceed algebraically in the quantum groups approach, we
work over $k$ a general field.

Next, also common, is the notion of `differential structure'
specified as an exterior algebra of differential forms $\Omega$
over $M$. An example is the universal differential calculus
$\Omega_{univ}$ going back to algebraic topology, Hochschild
cohomology etc. in the works of Quillen, Loday, Connes, Karoubi
and others. Here \eqn{Omegauniv}{\Omega^1_{univ}=\ker(\cdot),\quad
\extd f=f\tens 1-1\tens f,\quad \forall f\in M,} where
$\cdot:M\tens M\to M$ is the product. Similarly for higher
degrees. All other exterior algebras are quotients, so to classify
the abstract possibilities for $\Omega$ as a differential graded
algebra (DGA) we have to construct suitable differential graded
ideals. This is the line taken in \cite{Con} and an idea explained
there is to start with a spectral triple $(\rho(M),\CH,D)$ as a
representation of $\Omega_{univ}$, where $\extd$ is represented by
commutator with an operator $D$ called `Dirac', and then to divide
by the kernel. This does {\em not} however, usually, impose enough
constraints for a reasonable exterior algebra $\Omega(M)$, e.g. it
is often not finite-dimensional over $M$, so one forces higher
degrees to be zero, as well as additional relations with the aid
of a cyclic cocycle of the desired top degree $d$ corresponding to
a  `cycle' $\int :\Omega\to \C$. For an orientation there is also
a $\Z_2$ grading operator (trivial in the odd case) and a charge
conjugation operator, and Connes requires a certain Hochschild
cycle of degree $d$ whose image in the spectral triple
representation is the grading operator. This data together with
the operator norm and a K-theoretic condition for Poincar\'e
duality subsumes the role played in geometry by the usual Hodge *
operator and the volume form.

\begin{table}
\[\begin{array}{c|c|c}
{\rm Classical}& {\rm Connes\ approach}& {\rm Quantum\ groups\
approach}\\
\hline
 {\rm topol.\ space}& {\rm algebra}\ M& {\rm algebra}\ M\\
 {\rm integration}& \int:M\to \C, {\rm  trace} & \int:M\to k,
 \quad {\rm symmetry}\\
 {\rm differential\ calc.}& {\rm DGA}\ \Omega^\cdot,\ \extd & {\rm
 bimodule}\ \Omega^1\Rightarrow {\rm DGA}\ \Omega^\cdot\\
 {\ }\quad {\rm contruction\ of:}& {\rm spe.\ triple}\ (\rho(M),\CH,D)& {\rm
 classify\ all\ by\ symmetry}\\
 {\rm inner\ calculus}& -  & {\rm typ.}\ \exists \theta:\ \extd=[\theta, \}\\
 {\rm top\ form} & {\rm impose\ cycle}& {\rm typ.}\
 \exists\ {\rm Top}\\
\hline
{\rm principal\ bundle}& - & \begin{array}{c}(P,\Omega^1(P),A,\Omega^1(A),\Delta_R)\\
A\ {\rm quantum\ group\ fiber}\end{array}\\
{\rm vector\ bundle}& {\rm projective\ module}\
(\CE,e)& \CE=(P\tens V)^A,\ {\rm  typ.\ \exists\ }e\\
{\rm connection}& \nabla & \begin{array}{c}
\Omega^1(P)=\Omega^1_{hor}\oplus\Omega^1_{ver}\\
\Leftrightarrow\omega:\Lambda^1(A)\to \Omega^1(P) \Rightarrow
D_\omega\end{array}\\ {\rm Chern\ classes} &
{\rm Chern-Connes\ pairing}& -\\
\hline {\rm frame\ bundle}& - & {\rm princ.}\ P\ {\rm +\
soldering}\Rightarrow
\CE\isom\Omega^1(M)\\
{\rm metric}& ||[D,M ]||\Rightarrow d(\psi,\phi)&
g\in\Omega^1\tens_M\Omega^1;\ {\rm
coframing}\\
\begin{array}{c}{\rm Levi-Civta}\\{\rm Ricci\ tensor}
\end{array}&{\rm contained\ in}\  D&
\begin{array}{c}D_\omega\ +\ {\rm soldering}\Rightarrow \nabla_\omega\\
\extd\omega+\omega\wedge\omega\Rightarrow {\rm Ricci}\end{array}\\
\begin{array}{c}{\rm spin\ bundle}\\ {\rm Dirac\ operator}
\end{array}& {\rm assumed}\ (D,\CH)&\begin{array}{c} {\rm assoc.\ to\ frame\
bundle}\\  (\extd,\ \gamma_a,\omega)\Rightarrow \Dsl\end{array}\\
{\rm Hodge\ star}  & {\rm  in}\ ||\ ||\ +\ {\rm orientation} & g,\
{\rm Top}\ \Rightarrow *
\end{array}\]
\caption{Two approaches to noncommutative geometry}
\end{table}

Here again the quantum groups approach is less ambitious and we
don't have a general construction. Instead we {\em classify all
possible} $\Omega^1(M)$, typically restricted by some symmetry to
make the problem manageable. We postpone $\Omega^2(M)$ and higher
till later. Thus a first order calculus over $M$ is a nice notion
all by itself:
\begin{itemize}
\item An $M-M$-bimodule $\Omega^1$
\item A linear map $\extd:M \to \Omega^1$ such that
$\extd(fg)=(\extd f)g+f\extd g$ for all $f,g\in M$ and such that
the map $M\tens M\to \Omega^1$ given by $f\extd g$ is surjective.
\end{itemize}
It means classifying sub-bimodules of $\Omega^1_{univ}$ to
quotient by. In the last 5 years it was achieved on a case-by-case
bass for all main classes of quantum groups and bicovariant
calculi (see section~3). Then homogeneous spaces will likewise be
constrained by smoothness of desired actions and hence inherit
natural choices for their calculi, and so on for an entire
noncommutative universe of objects. Then we look at natural
choices for $\Omega^2(M)$, $\Omega^3(M)$  etc. layer by layer and
making choices only when needed.  There is a key
lemma\cite{BrzMa:dif} that if $\Omega^m(M)$ for $m\le n$ are
specified with $\extd$ obeying $\extd^2=0$ then these have a
maximal prolongation to higher degree, where we add in only those
relations at higher degree implied by $\extd^2=0$ and Leibniz. And
when $M$ is a quantum group there is a canonical extension to all
of $\Omega(M)$ due to \cite{Wor:dif} and the true meaning of which
is Poincar\'e duality in the sense of a non-degenerate pairing
between forms and skew-tensor fields. We will say more about this
in Section~3.

Apart from these differences in style the two approaches are again
broadly similar at the level of the differential forms and there
is much scope for convergence. Let us note one difference in
terminology. In the quantum groups approach we call a calculus
`inner' if there is a 1-form $\theta$ which generates $\extd$ as
graded commutator. This looks a lot like Connes `Dirac operator'
$D$ in the spectral triple representation \eqn{rhoD}{ \rho_D(\extd
f)=[D,\rho(f)]} but beware: our $\theta$ is a 1-form and nothing
to do with the geometric Dirac operator $\Dsl$ that we construct
only much later via gamma-matrices, spin connections etc. Moreover
$D$ in a spectral triple does {\em not} need to be in the image
under the representation $\rho_D$ of the space of 1-forms and
should not be thought of as a one form.

The next layers of geometry require bundles. In the traditional
approach to vector bundles coming from the theorems of Serre and
Swan one thinks of these as finitely generated projective modules
over $M$. However, in the quantum groups approach, as in
differential geometry, we think that important vector bundles
should really come as associated to principal ones. Principal
bundles should surely have a quantum group fiber so this is the
first place where the quantum group approach comes into its own.
The Brzezi\'nski-Majid theory of these is based on an algebra $P$
(the `total space'), a quantum group coordinate algebra $A$, {\em
and differential structures on them}, a coaction $\Delta_R:P\to
P\tens A$ such that $M=P^A$ and an exact sequence\cite{BrzMa:gau}
\eqn{bundle}{ 0\to P\Omega^1(M)P\to \Omega^1(P){\buildrel {\rm
ver}\over \longrightarrow} P\tens \Lambda^1(A)\to 0 } which
expresses `local triviality'. Here $\Lambda^1(A)$ are the
left-invariant differentials in $\Omega^1(A)$ and the map ver is
the generator of the vertical vector fields. We require its kernel
to be exactly the `horizontal' forms pulled up from the base. In
the geometrically less interesting case of the universal
differential calculus  the exactness of (\ref{bundle}) is
equivalent to a Hopf-Galois extension. A connection is a splitting
of this sequence and characterised by a connection form $\omega$.
Of course, we have associated vector bundles $\CE=(P\tens V)^A$
for every $A$-comodule $V$, analogous to the usual geometric
construction. A connection $\omega$ induces a covariant derivative
$D_\omega$ with the `derivation-like' property that one expects
for an abstract covariant derivative $\nabla$, so the two
approaches are compatible.

Just as the main example of the K-theory approach used to be the
`noncommutative torus' $T^2_\theta$\cite{Con:alg}, with vector
bundles classified in \cite{ConRie:yan}, the main example for
noncommutative bundles in the quantum groups approach was the
$q$-monopole bundle\cite{BrzMa:gau} over the Podle\'s
`noncommutative sphere' $S_q^2$.  Probably the first nontrivial
convergence between these approaches was in 1997 with the
Hajac-Majid projector $e$ where\cite{HajMa:pro}: \eqn{projSq}{
1-e=\left(\begin{matrix}q^2(1-b_3) & -qb_+ \cr b_-&
b_3\end{matrix}\right)\in M_2(S_q^2), \quad e \extd e\isom
D_\omega,\quad \<[\tau],[e]\>=\tau ({\rm Tr}(e))=1} where $\omega$
is the $q$-monopole connection and $\tau$ is the Masuda et al.
trace on $S_q^2$ which had been found in \cite{MMNNU:jac}. The
pairing is the Chern-Connes one between $HC^0$ and the class of
$e$ as an element of the K-group $K_0$. We use the standard
coordinates on $S_q^2$ namely with generators $b_\pm,b_3$ and
relations
\[ b_\pm b_3=q^{\pm 2}b_3 b_\pm+(1-q^{\pm 2})b_\pm,
\quad q^2b_- b_+=q^{-2}b_+ b_-+(q-q^{-1})(b_3-1)\] \eqn{Sqrel}{
b_3^2=b_3+q b_- b_+.} In this way one may exhibit the lowest
charge monopole bundle $\CE_1$ as a projective module. At the same
time, since the Chern-Connes pairing is nonzero, it means that the
$q$-monopole bundle is indeed nontrivial as expected. The
computations in \cite{HajMa:pro} are algebraic and, moreover, only
for the universal differential calculus, but we see that matching
up the two approaches is useful even at this level. This projector
has sparked quite a bit of interest in recent years. Dabrowski and
Landi looked for similar constructions and projectors for
$q$-instantons, replacing complex numbers by quaternions. An
inspired variant of that became the Connes-Landi projectors for
twisted spheres $S_\theta^4$. In these cases the new idea is (as
far as I understand) to work backwards from an ansatz for the form
of projector to the forced commutation relations for those matrix
entries as required by $e^2=e$ and a desired pairing with cyclic
cohomology. Note that, as twists, the noncommutative differential
geometry of such examples is governed at the algebraic level by
the Majid-Oeckl twisting theorem\cite{MaOec:twi}.

Finally, we come to the `top layer' which is Riemannian geometry.
In the approach of Connes this is all contained in the `Dirac
operator' $D$ which was assumed at the outset as an axiom. In the
quantum groups approach  we have been building up the different
layers and hope to construct a particular family of $\Dsl$
reflecting all the choices of differential structure, spin
connections etc. that we have made at lower levels. So the two
approaches are going in opposite directions. In the table we show
the quantum groups formulation of Riemannian geometry introduced
in \cite{Ma:rie} and studied recently for finite sets
\cite{Ma:rief}, based on the notion of a quantum frame bundle. The
main idea is a principal bundle $(P,A,\Delta_R)$ and an
$A$-comodule $V$ along with a `soldering form' $\theta_V:V\to
\Omega^1(P)$ that ensures that \eqn{frameisom}{ \Omega^1(M)\isom
\CE=(P\tens V)^A} i.e. that the cotangent bundle really can be
identified as an associated bundle to the principal `frame' one.
It turns out that usual notions proceed independently of the
choice of $A$, though what $\omega$ are possible and what
$\nabla_\omega$ they induce depends on the choice of $A$. Here
$\nabla_\omega:\Omega^1(M)\to \Omega^1(M)\tens_M\Omega^1(M)$ is
$D_\omega$ mapped over under the soldering isomorphism. A certain
$\bar D_\omega\wedge \theta_V$ corresponds similarly to its
torsion tensor. In terms of $\nabla_\omega$ the Riemann and
torsion tensor are given by \eqn{riem}{{\rm Riemann}=
(\id\wedge\nabla_\omega-\extd\tens\id)\nabla_\omega,\quad {\rm
Tor}=\extd - \nabla_\omega} as a 2-form valued operator an 1-forms
and a map from 1-forms to 2-forms, respectively. For the Ricci
tensor we need to lift the 2-form values of Riemann to
$\Omega^1(M)\tens_M\Omega^1(M)$ in a way that splits the
projection afforded by the wedge product, after which we can take
a trace over $\Omega^1$. This is possible at least in some cases.
This much is a `framed manifold' with connection.

A {\em framed Riemannian manifold} needs in addition a
nondegenerate tensor $g\in \Omega^1(M)\tens\Omega^1(M)$, which is
equivalent to a coframing $\theta_{V^*}:V^*\to \Omega^1(P)$ via
\eqn{gtheta}{ g=\<\theta_V\tens_M \theta_{V^*}\>,}where we pair by
evaluating on the canonical element of $V\tens V^*$,  and
$\theta_{V^*}$ is such that $\CE^*\isom \Omega^1(M)$. Since
$\CE^*$ under the original framing is isomorphic to
$\Omega^{-1}(M)$, such a coframing is the same as an isomorphism
$\Omega^{-1}(M)\isom\Omega^1(M)$. The coframing point of view is
nice because it points to a natural self-dual generalization of
Riemannian geometry: instead of demanding zero torsion and
$\nabla_\omega g=0$ as one might usually do for a Levi-Civita
connection, one can demand zero torsion and zero cotorsion (i.e.
torsion with respect to $\theta_{V^*}$ as framing). The latter is
a skew-version \eqn{skewcomp}{
(\nabla_\omega\wedge\id-\id\wedge\nabla_\omega) g=0} of metric
compatibility as explained in \cite{Ma:rie}. We can also demand
symmetry, if we want, in the form \eqn{gsymm}{ \wedge (g)=0} and
we can limit ourselves to $\theta_{V^*}$ built from $\theta_V$
induced by an $A$-invariant local metric $\eta\in V\tens V$. Let
us also note that if $M$ is parallelizable we can frame with a
trivial tensor product bundle and $\theta_{V},\theta_{V^*}$ reduce
to a vielbein $e_V:V\hookrightarrow \Omega^1(M)$ and covielbein
$e_{V^*}:V^*\hookrightarrow \Omega^1(M)$, i.e. subspaces forming
left and right bases respectively over $M$ and dual as
$A$-comodules. A connection $\omega$ now reduces to a `Lie
algebra-valued' 1-form \eqn{localcon}{ \alpha:\Lambda^1(A)\to
\Omega^1(M),} etc., in keeping with the local picture favoured by
physicists. One has to solve for zero torsion and zero cotorsion
in the form $\bar D\wedge e$ and $D\wedge e^*=0$. At the time of
writing the main noncommutative examples are when $M$ is itself a
quantum group. In the coquasitriangular case the dual of the space
of invariant 1-forms forms a braided-Lie algebra \cite{GomMa:bra},
which comes with a braided-Killing form $\eta$. This provides a
natural metric and in several examples one finds for it (by hand;
a general theorem is lacking) that there is a unique associated
generalised Levi-Civita connection in the sense
above\cite{Ma:rief}\cite{Ma:ric}.

We are then able to take a different $A$-comodule $W$, say, for
spinors. The associated bundle $\CS=(P\tens W)^A$ gets its induced
covariant derivative from the spin connection $\omega$ on the
principal bundle, and in many cases there is a reasonable choice
of `gamma-matrices' appropriate to the local metric $\eta$. We
then define the Dirac operator from these objects much as usual.
By now the approach is somewhat different from the Connes one and
we do {\em not} typically obtain something obeying the axioms for
$D$. This seems the case even for finite groups\cite{Ma:rief} as
well as for $q$-quantum groups\cite{Ma:ric}. The fundamental
reason is perhaps buried in the very notion of vector field: in
the parallelizable case an $M$-basis $\{e_a\}$ of $\Omega^1(M)$
implies `partial derivatives' $\del^a$ defined by $\extd f= \sum_a
e_a \del^a(f)$. These are not usually derivations but more
typically `braided derivations' (e.g. on a quantum group this is
shown in \cite{Ma:cla}). In cases such as the noncommutative torus
one has in fact ordinary derivations around. The noncommutative
differential calculus is a twist so that the constructions look
close to classical. But for $q$-examples and even finite group
examples, this is not at all the case. Perhaps this is at the root
of the mismatch and may stimulate a way to fix the problem.

Also in the presence of a metric we obtain a Hodge $*$ operator
$\Omega^m\to\Omega^{d-m}$ where $d$ is the `volume dimension' or
degree of the top (volume) form, assuming of course that it
exists. Once we have this we can write down actions such as
$-{1\over 4}F\wedge *F$ etc where $F=\extd \alpha$ is the
curvature of $\alpha\in\Omega^1$ modulo exact forms (Maxwell
theory) or $F=\extd \alpha+\alpha\wedge \alpha$ is the curvature
of $\alpha$ viewed modulo gauge transformation by invertible
functions (this is $U(1)$-Yang-Mills theory). Integration and the
Hodge * thus play the role of the operator norm $||\ ||$ in the
Connes approach. The two approaches were compared in a simple
model in \cite{MaSch:lat}. In principle one should be able to
extend these ideas to the non-Abelian gauge theory on bundles as
well, to construct a variety of Lagrangian-based models.

\section{Classifying calculi on general Hopf algebras}

We now focus for the rest of the paper on a small part of the
quantum groups approach discussed above, namely just the
differential calculus, and for the most part just $\Omega^1$. In
this section we let $M=A$ be a Hopf algebra over a field $k$.
Following Woronowicz\cite{Wor:dif}, a differential structure is
bicovariant if:

\begin{itemize}
\item $\Omega^1$ is a bicomodule with $\Delta_L:\Omega^1\to
A\tens \Omega^1$ and $\Delta_R:\Omega^1\to \Omega^1\tens A$
bimodule maps.
\item $\extd$ is a bicomodule map.

\end{itemize}

Here a Hopf algebra means a coproduct $\Delta:A\to A\tens A$, a
counit $\eps:A\to k$ and an antipode $S:A\to A$ such that $A$ is a
coalgebra, $\Delta$ an algebra map etc. \cite{Ma:book,Ma:prim}.
Coalgebras and (bi)comodules are defined in the same way as
algebras and (bi)modules but with the directions of structure maps
reversed. In the bicovariance condition $A$ is itself a
bi(co)module via the (co)product and $\Omega^1\tens A,
A\tens\Omega^1$ have the tensor product (bi)module structure.  The
second condition in particular fully determines
$\Delta_L,\Delta_R$ by compatibility with $\Delta$, so a
bicovariant calculus means precisely one where left and right
translation expressed by $\Delta$ extend consistently to
$\Omega^1$. The universal calculus $\Omega^1_{univ}\subset A\tens
A$ is bicovariant with coactions the tensor product of the regular
coactions defined by the coproduct on each copy of $A$.

The result in \cite{Wor:dif} is that $\Omega^1$ in the bicovariant
case is fully determined by the subspace $\Lambda^1$ of (say)
$\Delta_L$-invariant 1-forms. Indeed, there is a standard
bi(co)module isomorphism \eqn{hopfisom}{ A\tens A\isom A\tens
A,\quad a\tens b\mapsto a \Delta b} under which
$\Omega^1_{univ}\isom A\tens A^+$, where $A^+=\ker\eps$ is the
augmentation ideal or kernel of the counit (classically it would
be the functions vanishing at the group identity). The bimodule
structure on the right hand side of (\ref{hopfisom}) is left
multiplication in the first $A$ from the left and the tensor
product of two right multiplications from the right. The
bicomodule structure is the left coproduct on the first $A$ from
the left and the tensor product of right comultiplication and the
right quantum adjoint coaction from the right. Hence we arrive at
the classic result:

\begin{propos} (Woronowicz) Bicovariant $\Omega^1$ are in 1-1
correspondence with quotient objects $\Lambda^1$ of $A^+$ as an
$A$-crossed module under right multiplication and the right
quantum adjoint coaction.
\end{propos}

We recall that an $A$-crossed module means a vector space which is
both an $A$-module and a compatible $A$-comodule; the
compatibility conditions are due to Radford and correspond in the
finite-dimensional case to a module of the Drinfeld double
$D(A)=A^{*\rm op}\bowtie A$ when we view a right coaction of $A$
as a right action of $A^{*\rm op}$ by evaluation. Given the
crossed module $\Lambda^1$ we define $\Omega^1=A\tens\Lambda^1$
with the regular left(co)modules and the tensor product
(co)actions from the right. Because the category is prebraided,
there is a Yang-Baxter operator
$\Psi=\Psi_{\Lambda^1,\Lambda^1}:\Lambda^1\tens\Lambda^1\to\Lambda^1\tens\Lambda^1$
which is invertible when $A$ has bijective antipode. One has the
same results for right-invariant 1-forms $\bar\Lambda^1$ with
$\Omega=\bar\Lambda^1\tens A$.

Finally we recall that in the quantum groups approach we only need
to classify $\Omega^1$ because, at least for a bicovariant
calculus, there is a natural extension to an entire exterior
algebra $\Omega=A\tens \Lambda$ (where $\Lambda$ is the algebra of
left-invariant differential forms). There is a $\extd$ operation
obeying $\extd^2=0$ i.e. we have an entire DGA. The construction
in \cite{Wor:dif} is a quotient of the tensor algebra over $A$,
 \eqn{Worext}{
\Omega=T_A(\Omega^1)/\oplus_n \ker A_n,\quad A_n=\sum_{\sigma\in
S_n}(-1)^{l(\sigma)}\Psi_{i_1}\cdots\Psi_{i_{l(\sigma)}}} where
$\Psi_i\equiv\Psi_{i,i+1}$ denotes a certain braiding
$\Psi:\Omega^1\tens_A\Omega^1\to\Omega^1\tens_A\Omega^1$ acting in
the $i,i+1$ place of $(\Omega^1)^{\tens^n_A}$ and
$\sigma=s_{i_1}\cdots s_{i_{l(\sigma)}}$ is a reduced expression
in terms of simple reflections. This turns out to be equivalent to
defining $\Lambda$ directly as the tensor algebra of $\Lambda^1$
over $k$ with $A_n$ defined similarly by
$\Psi_{\Lambda^1,\Lambda^1}$ above. As such, $\Lambda$ is
manifestly a braided group or Hopf algebra in a braided category
of the linear braided space type\cite{Ma:book}. Similarly
$\Lambda^*$ starting from the tensor algebra of $\Lambda^{1*}$ is
a braided group, dually paired with $\Lambda$. The relations in
(\ref{Worext}) can be interpreted as the minimal such that this
pairing is nondegenerate, which is a version of Poincar\'e
duality. Also, $\Lambda$ being a braided group means among other
things that if there is a top form (an integral element) then it
is unique, which defines an `epsilon tensor' by the coefficient of
the top form. This combines with a metric to form a Hodge
* operator as mentioned in Section~2. More details and an application
to Schubert calculus on flag varieties are in our companion paper
in this volume.

\section{Classification on coquasitriangular Hopf algebras}

After Woronowicz's 1989 work for any Hopf algebra, the next
general classification result is my result for factorisable
coquasitriangular Hopf algebras such as appropriate versions of
the coordinate algebras $\C_q[G]$ of the Drinfeld-Jimbo quantum
groups, presented in Goslar, July 1996 \cite{Ma:gos,Ma:cla}. We
start with the finite-dimensional theory which is our main
interest for later sections, e.g. the above quantum groups at $q$
a primitive odd root of unity. We recall that a Hopf algebra $H$
is quasitriangular if there is a `universal R-matrix' $\CR\in
H\tens H$ obeying certain axioms (due to Drinfeld). It is
factorisable if $\CQ=\CR_{21}\CR$ is nondegenerate when viewed as
a linear map from the dual. The dual notion to quasitriangular is
coquasitriangular in the sense of a functional $\CR$ on the tensor
square, a notion which works well in the infinite-dimensional case
also\cite{Ma:book}.

\begin{theorem}(Majid) Let $A$ be a finite-dimensional factorisable
coquasitriangular Hopf algebra with dual $H$. Bicovariant
$\Omega^1(A)$ are in 1-1 correspondence with two-sided ideals of
$H^+$. \end{theorem}
\proof  This is  \cite[Prop. 4.2]{Ma:cla}
which states that the required action of the quantum double on
$\Lambda^{1*}\subset H^+$ under the quantum Killing form
isomorphism $\CQ:A^+\to H^+$ is the left and right coregular
representation of two copies of $H$ on $A^+$ (the coregular action
of $H$ on $A$ is obviously adjoint to the regular action of $H$ on
$H$ if we wish to phrase it explicitly in terms of $H$, i.e.
$\Lambda^1$ is isomorphic to a quotient of $H^+$ by a 2-sided
ideal). \eproof

As far as the actual proof in \cite{Ma:cla} is concerned, there
were two ideas. First of all, instead of classifying quotients
$\Lambda^1$ of $A^+$ it is convenient to classify their duals or
`quantum tangent spaces' as subcrossed modules $\CL\subset H^+$.
For coirreducible calculi (ones with no proper quotients) we want
irreducible $\CL$, but we can also classify indecomposable ones,
etc. Of course an $A$-crossed module can be formulated as an
$H$-crossed module (the roles of action and coaction are swapped)
or $D(H)$-module, so actually we arrive at a self-contained
classification for quantum tangent spaces for any Hopf algebra
$H$. One has to dualise back (as above) to get back to the
left-invariant 1-forms. The second idea was that when $H$ is
factorisable the quantum Killing form is a nondegenerate map
$\CQ:A^+\to H^+$ and, moreover, $D(H)\isom H\codcross H$, which as
an algebra is a tensor product $H\tens H$. We refer to
\cite{Ma:book} where the forward direction of the isomorphism was
proven for the first time. So the crossed module structures
$\CL\subset H^+$ that we must classify become submodules of $A^+$
under the action of $H$ from the left and the right (viewed as
left via $S$), which we computed as the left and right coregular
ones. Converting back to $\Lambda^1$ means of course 2-sided
ideals as stated.

At present we are interested in the finite-dimensional case where
we need only the algebraic theorem as above. Note that every
Artinian algebra has a unique block decomposition, which includes
all finite-dimensional algebras $H$ over a field $k$. The
decomposition is equivalent to finding a set of orthogonal
centrally primitive idempotents $e_i$ with \eqn{projsum}{ 1=\sum_i
e_i.} These generate ideals $e_iH$. Note that $e_i^2=e_i$ implies
that $\eps(e_i)=0,1$ and the above implies that exactly one is
nonzero. Similarly $H^+=\oplus_i e_i H^+$ is a decomposition of
$H^+$. Hence \eqn{univdec}{ \Omega^1_{univ}\isom\oplus_i
\Omega^1_{e_i}} where, for any central projector $e$ we have a
calculus \eqn{Omegae}{ \Omega^1_{e}= A\tens\Lambda^1_{e},\quad
\Lambda^1_{e}=eH^+.} We build the left-invariant 1-forms directly
on the block as isomorphic to $H^+$ modulo the kernel of
multiplication by $e$. In these terms (tracing through the details
of the classification theorem\cite{Ma:cla}) we have explicitly:
\eqn{ecalc}{ (eh).a=\sum a\o e \CR_2(a\t)h\CR_1(a\thr),\quad \extd
a=\sum a\o e\CQ_1(a\t)-a e;\quad \theta=e} (the calculus here is
inner). We use the Sweedler notation $\Delta a=\sum a\o\tens a\t$
and the notation from \cite{Ma:book} where $\CR_1(a)=(a\tens
\id)(\CR)$ etc. Since $\CR_1$ is an algebra map and $\CR_2$ an
antialgebra map, one may easily verify that this defines a
bimodule and that the Leibniz rule is obeyed. These formulae
generalise ones given usually in terms of R-matrices. One has, cf.
\cite{Ma:cla}:

\begin{corol}(Majid) In the finite-dimensional semisimple case, coirreducible
bicovariant calculi on $A$ in the setting of Theorem~4.1 are in
1-1 correspondence with nontrivial irreducible representations of
the dual Hopf algebra $H$. The dimension of the calculus is the
square of the dimension of the corresponding representation.
\end{corol}
\proof Indeed, in the case of $H$ finite-dimensional semisimple,
each block will be a matrix block corresponding to an irreducible
representation. In this case among the projectors there will be
exactly one where $eh=e\eps(h)$ for all $h\in H$ (the normalized
unimodular integral). It has counit 1 and corresponds to the
trivial representation; we exclude it in view of $eH^+=0$. \eproof

This semsimple case is the situation covered in \cite{Ma:cla} in
the form of an assumed Peter-Weyl decomposition. We will apply it
to the quantum double example in the next section. Moreover, for
each irredicble representation, one may write (\ref{ecalc}) in
terms of the matrix $R$ given by $\CR$ in the representation. In
this case one has formulae first used by Jurco\cite{Jur} for the
construction of bicovariant calculi on the standard quantum groups
such as $\C_q[SL_n]$. That is not our context at the moment but we
make some remarks about it at the end of the section. Let us note
only that at the time they appeared such results in \cite{Ma:cla}
were the first of any kind to identify the full moduli of all
coirreducible bicovariant calculi in some setting with irreducible
representations.

In the nonsemisimple case the algebra $H$ has a nontrivial
Jacobson radical $J$ defined as the intersection of all its
maximal 2-sided ideals. It lies in $H^+$. Hence by Theorem~4.1
there is a calculus \eqn{Omss}{ \Omega_{ss}^1=A\tens H^+/J} which
we call the `semisimple quotient' of the universal calculus. $H/J$
has a decomposition into matrix blocks giving a decomposition of
$\Omega^1_{ss}$ along the lines of the semisimple case.

\begin{example} Let $\C^{\rm red}_q[SL_2]$ be the 27-dimensional reduced quantum group
at $q$ a primitive cube root of unity. Here
\[ a^3=d^3=1,\quad b^3=c^3=0\]
in terms of the usual generators. The enveloping Hopf algebra
$u_q(sl_2)$ is known to have the block decomposition
\[ u_q(sl_2)=M_3(\C)\oplus B\]
where $B$ is an  18-dimensional non-matrix block (the algebra is
not semisimple) with central projection of counit 1. Hence the
universal calculus decomposes into nontrivial calculi of
dimensions $9,17$. Moreover, $J\subset B$ is 13-dimensional and
the quotient
\[ B/J=\C\oplus M_2(\C)\]
implies a 4-dimensional calculus as a quotient of the
17-dimensional one (the other summand $\C$ gives zero). Here
$\Omega^1_{ss}$ is the direct sum of the 4 and 9 dimensional
matrix calculi.\end{example}

The natural choice of calculus here is 4-dimensional and has the
same form as the lowest dimension calculus for generic $q$, the 4D
one first found by hand by Woronowicz\cite{Wor:dif}. For roots of
unity the cohomology and entire geometry are, however, completely
different from the generic or real $q$ case as shown in
\cite{GomMa:coh}. This takes $\C^{\rm red}_q[SL_2]$ at 3,5,7-th roots as a
finite geometry where all computations can be done and all ideas
explored completely. We find, for example, the Hodge * operator
for the natural $q$-metric, and show that the moduli space of
solutions of Maxwell's equations without sources or
`self-propagating electromagnetic modes' decomposes into a direct
sum
\[ \{{\rm Maxwell\ zero\ modes}\}=\{ {\rm zero-curvature}\}
\oplus\{{\rm self-dual}\}\oplus\{{\rm anti-selfdual}\}\] of
topological (cohomology) modes, self-dual and anti-selfdual ones.
The noncommutative de Rham cohomology here in each degree has the
same dimension as the space of left-invariant forms, for reasons
that are mysterious. We also find that the number of self-dual
plus zero curvature modes appears to coincide with the number of
harmonic 1-form modes. We find in general that the reduced
$\C^{\rm red}_q[SL_2]$, although finite dimensional and totally
algebraic, behaves geometrically like a `noncompact' manifold
apparently linked to the nonsemisimplicity of $u_q(sl_2)$. Also
interesting is the Riemannian geometry of the reduced $\C^{\rm
red}_q[SL_2]$ in \cite{Ma:ric}. The Ricci tensor for the same
$q$-metric turns out to be essentially proportional to the metric
itself, i.e. an `Einstein space'. It would also be interesting to
look at the differential geometry induced by taking the other
blocks according to our classification theorem, particularly the
non-matrix block.

Finally, let us go back and comment on the situation for the usual
(not reduced) $q$-deformation quantum groups $\C_q[G]$ associated
to simple Lie algebras $\cg$. There are two versions of interest,
Drinfeld's deformation-theoretic formal powerseries setting over
$\C[[\hbar]]$ where $q=e^{\hbar\over 2}$, and the algebraic fixed
$q$ setting. The first is easier and we can adapt the proof of
Theorem~4.1 immediately to it\cite{Ma:cla}. Here $H=U_q(\cg)$ and
$A=U_q(\cg)'$ is its topological dual Hopf algebra rather than the
more usual coordinate algebra (the main difference is that some of
its generators have logarithms).

\begin{propos}(Majid) In the formal powerseries setting over $\C[[\hbar]]$,
bicovariant $\Omega^1$ on the dual of $U_q(\cg)$ are in 1-1
correspondence with 2-sided ideals in $U(\cg)[[\hbar]]^+$.
\end{propos}
\proof  In fact the key point of Drinfeld's formulation \cite{Dri}
is that one has a category of deformation Hopf algebras with all
the axioms of usual Hopf algebra including duality holding over
the ring $\C[[\hbar]]$ and the additional axioms of a
quasitriangular structure holding in particular for his version of
$U_q(\cg)$. The arguments in the proof of \cite[Prop. 4.2]{Ma:cla}
and Theorem~4.1, as outlined above, require only these general
axiomatic properties and take the same form line by line in the
Drinfeld setting. Thus, the map $\CQ$ has the same intertwiner
properties used to convert an ad-stable right ideal in the dual to
a 2-sided ideal in $U_q(\cg)$. That $U_q(\cg)$ are indeed
factorizable (so that $\CQ$ is an isomorphism) when we use the
topological dual is the Reshetikhin-Semenov-Tian-Shanksy
theorem\cite{ResSem:mat}. Its underlying reason is visible at the
level of the subquantum groups $U_q(b_\pm)$ where
$U_q(b_\pm)'\isom U_q(b_\mp)$ was shown by Drinfeld in \cite{Dri},
which combines with the triangular decomposition of $U_q(\cg)$ to
yield factorisibility. Finally, also due to Drinfeld, is that
$U_q(\cg)\isom U(\cg)[[\hbar]]$ as algebras (since given by a
coproduct-twist of the latter as a quasiHopf
algebra\cite{Dri:qua}). Hence two sided ideals in $U_q(\cg)^+$ are
in correspondence with 2-sided ideals in the undeformed
$U(\cg)[[\hbar]]^+$. \eproof

This reduces the classification of calculi to a classical question
about the undeformed algebra $U(\cg)[[\hbar]]$. The natural
two-sided ideals of interest here are those given by two-ideals
from $U(\cg)$ as a Hopf algebra over $\C$. Let us call these
ideals of `classical type' and likewise the corresponding calculi
`of classical type'. We are interested in those contained in
$U(\cg)^+$ and which are cofinite so that the corresponding
calculi are finite-dimensional (f.d.). Maximal such ideals
correspond to coirreducible calculi. Maximal cofinite ideals of
$U(\cg)$, as for any algebra over a field, correspond (via the
kernel) to its finite-dimensional irreducible representations and
hence to such representations of $\cg$ as a Lie algebra.
Intersection with $U(\cg)^+$ gives corresponding ideals which are
maximal there, and for a proper ideal we drop the trivial
representation. The correspondence here uses the central character
to separate maximal ideals in one direction, and the annihilator
of the quotient of $U(\cg)^+$ by the ideal in the other direction.
Hence Proposition~4.4 has the corollary \eqn{calcg}{ \left\{{{\rm
coirreducible\ classical-type}\atop{\rm f.d.\ bicovariant\
calculi\ on}\ U_q(\cg)'}\right\}\leftrightarrow \left\{{{\rm
nontrivial\ f.d.\ irreducible}\atop{\rm representations\ of}\
\cg}\right\}} which is in the same spirit as Corollary~4.2. Let us
note that when one speaks of representations of $U_q(\cg)$ in the
deformation setting one usually has in mind ones similarly of
classical-type deforming ones of $U(\cg)$ (so that one speaks of
an equivalence of categories). On the other hand, without a
detailed study of the Peter-Weyl decomposition for $U_q(\cg)$ in
the $\C[[\hbar]]$ setting, we do not claim that the universal
calculus decomposes into a direct sum of these calculi, which
would be the full-strength version of the theory in \cite{Ma:cla}
(as explained there).

The classification for the $q\in \C$ case with $\C_q[G]$ algebraic
is more complicated and shows a different `orthogonal' kind of
phenomenon. Here one finds, as well as one calculus for each
finite-dimensional irreducible representation (in the spirit of
(\ref{calcg})), further `twisted variants' of the same square
dimension. This possibility was mentioned in \cite{Ma:cla} and
attributed to the fact that in this case the Hopf algebra is not
quite factorisable. Points of view differ on the significance of
these additional twists but our own is the following. First of
all, one can guess that the calculi that fit with the
$\C[[\hbar]]$ analysis above should correspond in the algebraic
setting to calculi which are commutative as $q\to 1$, while the
`twisted variants' should not have this property and hence could
be considered as pathological from a deformation-theoretic point
of view. As evidence, we demonstrate this below for $\C_q[SL_n]$.
Moreover, this is a pathology that exists for most R-matrix
constructions, not only calculi. Indeed, we already explained in
\cite{Ma:qua,Ma:book} that to fit with Hopf algebra theory
R-matrices must be normalised in the `quantum group
normalisation'. For $\C_q[SL_n]$ let $R_{hecke}$ be in the usual
Hecke normalisation where the braiding has eigenvalues
$q,-q^{-1}$. Then the correct quantum group normalisation is
\eqn{R}{ R=q^{-{1\over n}}R_{hecke}} and we explained that one has
the ambiguity of the choice of $n$-th root. The Jurco construction
for the $n^2$-squared calculus was given entirely in terms of
R-matrices (we give some explicit formulae in Section~5) so one
has an $n$-fold ambiguity for the choice of $q^{-{1\over n}}$. The
unique choice compatible with the $\C[[\hbar]]$ point of view is
the principal root that tends to 1 as $q\to 1$. The other choices
differ by $n$-th root of unity factors in the normalisation and
thereby in all formulae and these are the `additional twists'.
Similarly in the contragradient and other representations. This
reproduces exactly what was found for the specific analysis of
calculi on $\C_q[SL_n]$ of dimension $n^2$ in \cite{SchSch:cla}
(note, however, that this was not a classification of all calculi
in the sense above since the dimension was fixed at $n^2$ or
less). It was found that for $n>2$ there were $2n$ calculi
labelled by $\pm$ and $z$ a primitive $n$th root of unity (or just
the parameter $z$ if $n=2$). As remarked already in
\cite{SchSch:cla}, only the $z=1$ cases have a commutative limit
as $q\to 1$ and we identify them now as the canonical choices for
the fundamental and conjugate fundamental representations (these
being identified if $n=2$). One can also view this at the level of
the universal R-matrix for all representations as a variation of
$\CR$.

Following \cite{Ma:cla} (their preprint appeared several months
after \cite{Ma:cla} was archived), Bauman and Schmidt\cite{BauSch}
studied the classification for $\C_q[G]$ using the same
factorizablity method as above. This is a more formal treatment
but not a complete analysis of the possible twists any more than
\cite{Ma:cla} was. Meanwhile, about a year after \cite{Ma:cla},
Heckenberger and Schm\"udgen\cite{HecSch:cla} gave a full
classification for certain but not all $\C_q[G]$ using a different
method, which appears to be the current state of play.

\section{Differentials on the quantum codouble of a finite group.}

In this section we are going to demonstrate the classification
theorem in the previous section for the most famous factorisable
coquasitriangular Hopf algebra of all, namely the coordinate
algebra of the Drinfeld quantum double itself. So
$A=D^*(\CA)=\CH^{cop}\codcross \CA$ where $\CA$ is a finite
dimensional Hopf algebra and $\CH$ is its dual. As an algebra the
codouble is a tensor product but we are interested in bicovariant
calculi on it, which depends on the doubly-twisted coproduct. By
our theorem, these are classified by two-sided ideals in $D(\CH)$.
Actually, we compute only the case where $\CA=k(G)$, the functions
on a finite group $G$, but the methods apply more generally. Then
$A=k G\lcocross k(G)$. We assume $k$ is of characteristic zero.

\begin{theorem} Differential calculi on $D^*(G)=kG \lcocross k(G)$ are
classified by pairs $(\CC,V)$ where $\CC\subset G$ is a conjugacy
class and $V$ is an irreducible representation of the centralizer,
and at least one of $\CC,V$ are nontrivial. The calculus has
dimension $|\CC |^2\dim(V)^2$.
\end{theorem}
\proof Here $H=D(G)=k(G)\lcross kG$ is a semidirect product by the
adjoint action. Taking basis $\{\delta_s\tens u|\ s,u\in G\}$, the
product is $(\delta_s\tens u)(\delta_t\tens
v)=\delta_{s,utu^{-1}}(\delta_s\tens uv)$. Consider an element of
the form $e=\sum_s \delta_s \tens e_s$ where $e_s\in kG$. To be
central, we consider
\[ \delta_t.e=\delta_t\tens
e_t,\quad e.\delta_t=\sum_s\delta_s\delta_{Ad_{e_s}(t)}\tens
e_s.\] Writing $e_s=\sum_u e_{s,u}u$, say, equality requires
\[ e_{t,u}\delta_t=e_{utu^{-1},u}\delta_{utu^{-1}},\quad \forall
t, u\in G.\] This implies that $e_s\in k G_s$ the group algebra of
the centralizer of $s$ in $G$. Next, $u.e=e.u$ requires that
$e_{usu^{-1}}=u e_s u^{-1}$ for all $u,s\in G$. Thus central
elements $e$ are of the form \eqn{eD}{ e=\sum_{s\in \CC}
\delta_s\tens e_s,\quad e_s\in kG_s, \quad
ue_su^{-1}=e_{usu^{-1}},\quad \forall u\in G} for some Ad-stable
subset $\CC\subset G$. In that case one may compute
$e^2=\sum_s\delta_s\tens e_s^2$. Hence projectors are precisely of
the above form with each $e_s$ a central idempotent of $kG_s$. For
a centrally primitive idempotent we need $\CC$ a conjugacy class
and a centrally primitive idempotent $e_0$ on the group algebra of
the centralizer $G_{0}$ (any one point $s_0\in \CC$ determines the
rest). The choice of $e_0$ comes from the block decomposition of
$kG_0$ and in characteristic zero this is given by irreducible
representations $V$. The relation is \eqn{e0}{ e_0={\dim(V)\over
|G_0|}\sum_{u\in G_0}\trace_V(u^{-1})u.} This gives the block
decomposition of $D(G)$. For $D(G)^+$ we remove the trivial case.
\eproof

Note that this result is exactly in line with Corollary 4.2 since
the result is the same data as for the classification of
irreducible representations of $D(G)$, which are also the same as
irreducible crossed $G$-modules. This is to be expected as $D(G)$
is semisimple because the square of its antipode is the identity.
Let $\{e_i\}$ be a basis of $V$ and $\CC=\{a,b,c,\dots\}$. Note
also that $\CC$ is just the data for a calculus $\Omega^1(k(G))$
on a finite group function algebra, while $V$ is just the data for
a calculus $\Omega^1(kG_0)$ on a finite group algebra viewed `up
side down' as a noncommutative space, a result in \cite{Ma:cla}.
The calculus on the double glues together these calculi. We now
obtain explicit formulae as follows.

Given the data above, the associated representation $W$ of $D(G)$
has basis $\{e_{ai}\}$ say, with action \[
\delta_s.e_{ai}=\delta_{s,a}e_{ai},\quad u.e_{ai}=\sum_j
e_{uau{-1}j}\, \zeta_a(u)^j{}_i,\] where $\zeta:\CC\times G\to
G_0$ is a cocycle \eqn{zeta}{\zeta_a(u)=g^{-1}_{uau^{-1}}u
g_a;\quad \zeta_a(uv)=\zeta_{vav^{-1}}(u)\zeta_a(v)} defined by
any section map $g:\CC\to G$ such that $g_a s_0 g_a^{-1}=a$ for
all $a\in \CC$,  and we use its matrix in the representation $V$.
Completely in terms of matrices, we have the representation of
$D(G)$ by \eqn{rhoDmat}{ \rho(\delta_s\tens
u)^{ai}{}_{bj}=\delta_{s,a}\delta_{u^{-1}au,b}\zeta_{b}(u)^i{}_j.}
We also need the universal R-matrix and quantum Killing form
\eqn{univD}{\CR=\sum_{u\in G}\delta_u\tens 1\tens 1\tens u,\quad
\CQ=\sum_{u,v}\delta_{uvu^{-1}}\tens u\tens \delta_u\tens v} of
$D(G)$. Finally, we need the Hopf algebra structure of
$A=kG\lcocross k(G)=\{s\tens\delta_u|\ s,u\in G\}$ with coproduct
\eqn{D*}{ \Delta \delta_u=\sum_{vw=u}\delta_v\tens\delta_w,\quad
\Delta s=\sum_u s\delta_u\tens u^{-1}su} and $s,\delta_u$
commuting.

We insert these formulae into the Jurco-type construction for a
representation $(\rho,W)$ and the conventions of \cite{Ma:cla}.
One can do it either for left-invariant forms $\Lambda$ as
elsewhere in the paper, or, which gives nicer results for our
present conventions for $D^*(G)$, we can do it for right invariant
forms $\bar\Lambda^1$. In this case we define
$\bar\Lambda^1=\End(W)$ with basis $\{e_{\alpha}{}^{\beta}\}$ and
use the standard formulae for quasitriangular Hopf algebras, cf.
Section~4:
\[\Omega^1=\End(W)\tens A,\quad a.e_{\alpha}{}^{\beta}=\sum
\rho(\CR_1(a\o))^{\gamma}{}_{\alpha}\, e{}_{\gamma}{}^{\delta} \,
\rho(\CR_2(a\t)){}^{\beta}{}_{\delta}\, a\thr\] \eqn{JurR}{\extd
a=\sum \rho(\CQ_2(a\o))a\t-\theta a,\quad \forall a\in A;\quad
\theta=\sum_{\alpha}e_{\alpha}{}^{\alpha},} where
$\rho(h)=\rho(h)^\alpha{}_\beta e_\alpha{}^\beta$ and summations
of indices are understood. Putting in the explicit formulae
(\ref{rhoDmat})--(\ref{D*}) we have basis $\{e{}_{ai}{}^{bj}\}$
say, and find \[ \Omega^1=\End(W).D^*(G),\quad f
e_{ai}{}^{bj}=e_{ai}{}^{bj} L_b(f),\quad \extd f= \sum_{a,i}
e_{ai}{}^{ai} \del_a(f)\] \eqn{OmD}{s e_{ai}{}^{bj}=\sum_k
e_{sas^{-1}k}{}^{bj}\, \zeta_a(s)^k{}_i b^{-1}sb,\quad \extd
s=\sum_{a,i,j} e_{sas^{-1}i}{}^{aj}\, \zeta_a(s){}^i{}_j
a^{-1}sa-\theta s} for all $s\in G, f\in k(G)$. Here
$\del_a=L_a-\id$,  where $L_a(f)=f(a(\ ))$ denotes
left-translation in the direction of $a$ and
$\theta=\sum_{ai}e_{ai}{}^{ai}$ makes the calculus inner. The
special case where we take the trivial representation of the
centralizer is a canonical calculus of dimension $|\CC|^2$
associated to any conjugacy class, cf. some similar formulae from
a different point of view (viewing the double as a bicrossproduct)
in \cite{NME:cla}.

\begin{example} For $S_3$ with $u=(12),v=(23),w=(13)$, the possible conjugacy
classes are  $\{e\},\{uv,vu\},\{u,v,w\}$ of orders 1,2,3. Their
centralizers are $S_3,\Z_3,\Z_2$. Hence on $D^*(S_3)$ we have
calculi of dimension 1,4 for the first class (the nontrivial
irreducibles of $S_3$), three calculi of dimension 4 (the three
irreducibles of $\Z_3$) for the second class and two of dimension
9 (the two irreducibles of $\Z_2$) for the third class.
\end{example}

Of these, we now focus on the 9-dimensional calculus since the
associated conjugacy class of order 3 defines the usual
differential calculus on $S_3$ and the $D^*(S_3)$ calculus is an
extension of that. As basepoint we take $s_0=u$ and as section we
take $g_u=e$, $g_v=w$, $g_w=v$. We let $q=\pm1$ according to the
trivial or nontrivial representation of $\Z_2$. The resulting
cocycle as a function on
$S_3$ is \eqn{DS3z}{\begin{array}{c|cccccc} S_3& 1&u&v&w&uv&vu\\
\hline
\zeta_u&1&q&1&1&q&q\\ \zeta_v  &1 &q &q &1 &1&q \\
\zeta_w  & 1&q&1&q&q&1\end{array}}

\medskip\noindent For the calculus, we obtain from the above that
$\extd f=\sum_a e_a\del_a(f)$ is the usual calculus on $S_3$ with
$e_a\equiv e_a{}^a$, and \eqn{DS3d}{ \extd u=q(e_u u+e_w{}^v
w+e_v{}^w v)-\theta u,\quad \extd v=e_u{}^wu+q e_v v+e_w{}^u
w-\theta v} and the same with $v,w$ interchanged in the second
expression. Here $\theta=e_u+e_v+e_w$. We also have
\begin{eqnarray*}&& \extd(uv)=(qe_u{}^w+qe_v{}^u+e_w{}^v)vu-\theta
uv\\
&& \extd(vu)=(qe_u{}^v+e_v{}^w+qe_w{}^u)uv-\theta
vu.\end{eqnarray*} Of course, it is enough to work with generators
$u,v$ say of $S_3$ using the commutation relations from
(\ref{OmD}), namely \eqn{DS3rel}{ ue_a{}^b=q e_{uau}{}^b bub,\quad
ve_u{}^b=e_w{}^b bvb,\quad ve_v{}^b=q e_v{}^b bvb,\quad
ve_w{}^b=e_u{}^b bvb.}

Finally, for physics, we also need the higher exterior algebra as
explained in Section~3. For coquasitriangular Hopf algebras and in
the present right-invariant conventions we have the standard
braiding \eqn{braidJur}{ \Psi(e_{\alpha}{}^{\beta}\tens
e_{\gamma}{}^{\delta})=e_\mu{}^\nu\tens
e_{\sigma}{}^{\tau}(R^{-1})^{\alpha_1}{}_\alpha{}^\mu{}_{\alpha_2}
R^\beta{}_{\alpha_3}{}^{\alpha_2}{}_\gamma
R^\delta{}_{\alpha_4}{}^{\alpha_3}{}_\tau
\tilde{R}^{\alpha_4}{}_\nu{}^\sigma{}_{\alpha_1}} where
$R=(\rho\tens\rho)(\CR)$ and $\tilde{R}=(\rho\tens\rho\circ
S)(\CR)$. This is adjoint to the braided-matrix relations of the
braided matrices $B(R)$, i.e.  $\bar\Lambda^{1*}$ is a standard
matrix braided Lie algebra\cite{GomMa:bra}. We can compute this
braiding explicitly for $D^*(G)$ above to find:
\begin{eqnarray}\label{braidD} &&\kern -33pt \Psi(e_{ai}{}^{bj}\tens
e_{ck}{}^{dl})\nonumber\\
&&=e_{a^{-1}bcb^{-1}a m}{}^{dl} \zeta_c(a^{-1}b)^m{}_k \tens
\zeta^{-1}_b(d^{-1})^j{}_p e_{d^{-1}ad n}{}^{d^{-1}bd
p}\zeta_a(d^{-1})^n{}_i\end{eqnarray} where we sum over the matrix
indices $m,n,p$ of the representation of the centralizer group.
The relations of the exterior algebra are then computed using
braided-factorial matrices \cite{Ma:book}.

Then for our example $D^*(S_3)$ with the standard order 3
conjugacy class and either representation $q=\pm 1$ of the
isotropy group $\Z_2$, we have dimensions and cohomology in low
degree:
\[ \dim(\Omega(D^*(S_3)))=1:9:48:198:\cdots, \]
\eqn{DS3coh}{\quad H^0(D^*(S_3)) =k.1,\quad
H^1(D^*(S_3))=k.\theta} which is the same as for an isomorphic
bicrossproduct example computed in \cite{NME:cla}. We expect the
exterior algebra to be finite-dimensional and the Hilbert series
to have a symmetric form (our computer did not have enough memory
for further degrees to check this). The exterior algebra appears
to be quadratic with the relations in degree 2 from (\ref{braidD})
as follows. For simplicity we state them only for $q=1$:
\[ e_a{}^b\wedge e_{aba^{-1}}{}^b=0,\quad (e_{aba^{-1}}{}^a)^2
+\{e_a{}^a,e{}_b{}^a\}=0,\quad \forall a,b\] \eqn{DS3exta}{
e_u\wedge e_v+e_v\wedge e_w+e_w\wedge e_u=0} and the conjugate
(product-reversal) of this. Here $e_a\equiv e_a{}^a$ obey
$e_a^2=0$ and we recover precisely the usual $\Omega(S_3)$ as a
subalgebra generated by them. In addition, we have
\[ e_u{}^v\wedge e_v{}^u+e_v{}^w\wedge e_w{}^v+e_w{}^u\wedge
e_u{}^w=0\] \eqn{DS3extb}{e_u{}^v\wedge e_u{}^w+e_v{}^w\wedge
e_v{}^u+e_w{}^u\wedge e_w{}^v=0} and their opposites, and
\eqn{DS3extc}{ e_u\wedge e_u{}^w+e_u{}^w\wedge e_v+e_u{}^v\wedge
e_w{}^u=0} and all permutations of $u,v,w$ in this equation, plus
all their conjugates.

\section{Calculus on general $D^*(\CA)$ and T-duality}

Here we describe the differential calculus on a general quantum
codouble Hopf algebra $D^*(\CA)$ associated to any $\CA$-crossed
module. These are not in general all calculi (they are the block
decomposition of the semisimple calculus $\Omega_{ss}$) and we do
not attempt a classification as we did for $D^*(G)$ above. Let
$\CA$ for the moment be a finite-dimensional Hopf algebra with
dual $\CH$ and $A=D^*(\CA)=\CH^{cop}\codcross \CA$ as in
Section~5. Its coproduct and coquasitriangular structure are, in
the conventions of \cite{Ma:book},
\[ \Delta(h\tens a)=\sum h\o\tens f^a a\o f^b\tens (Se_a)h\o e_b\tens
a\t\]
\[ \CR(h\tens a,g\tens b)=\eps(a)\eps(g)\<h,b\>,\quad \forall
h,g\in \CH,\quad a,b\in\CA\] where $\{e_a\}$ is a basis of $\CH$
and $\{f^a\}$ a dual basis. Meanwhile,  $D(\CH)=\CA^{op}\bowtie
\CH$ has the product
\[ (a\tens h)(b\tens g)=\sum b\o a\tens h\t
g\<Sh\o,b\o\>\<h\thr,b\thr\>\] as in \cite{Ma:book}. According to
Corollary~4.2 we get a matrix block calculus for any
representation of $D(\CH)$, which means a (left) $\CH$-crossed
module $W$. We let $\{e_\alpha\}$ be a basis for $W$,
$\{f^\alpha\}$ a dual basis and $e_\alpha{}^\beta=e_\alpha\tens
f^\beta$ as in the previous section. We write the left crossed
module action as $\la$ and coaction as $\sum w\bo\tens w\bt$.

\begin{propos} Let $W$ be an $\CH$-crossed module.  The corresponding
calculus $\Omega^1(D^*(\CA))=\End(W)\tens D^*(\CA)$ has the
relations
\begin{eqnarray*} &\extd h=\sum  (h\t\la e_\beta\bt\tens f^\beta)
\Ad_{e_\beta\bo}(h\o)-\theta h\\
 &\extd a=\sum \<\Ad_{f^a}(a\o),e_\beta\bo\>(e_\beta\bt\tens
f^\beta)e_a\, a\t  -\theta a\\
& h.e_\alpha{}^\beta=\sum h\t\la e_\alpha{}\tens f^\gamma
\Ad_{e_\gamma\bo}(h\o)\<e_\gamma\bt,f^\beta\>\\
 & a.e_\alpha{}^\beta=\sum e_\alpha{}\tens f^\gamma
\<\Ad_{f^a}(a\o),e_\gamma\bo\>e_a \<e_\gamma\bt,f^\beta\> a\t
\end{eqnarray*} for all $a\in \CA$ and $h\in\CH$, where $\Ad$ is the right
adjoint action and $\theta=\sum e_\alpha{}^\alpha$.
\end{propos}
\proof  From the definition of the matrices $\rho$ we have
\[ \rho(a)=\sum \<a,e_\alpha\bo\>e_\alpha\bt\tens f^\alpha,\quad \rho(h)=\sum h\la
e_\alpha\tens f^\alpha.\] We then use the above formulae for the
quantum double and find in particular the required map
\[ \CQ_2(h\tens a)=(1\tens h)(a\tens 1)=\sum a\t\tens h\t
\<Sh\o,a\o\>\<h\thr,a\thr\>.\] We then use (\ref{JurR}), making
routine Hopf algebra computations, including (for $\extd h$) the
crossed-module compatibility conditions \[ \sum (h\o\la g)\bo
h\t\tens (h\o\la g)\bt=\sum h\o g\bo\tens h\t\la g\bt\] to obtain
the result. We use $\Ad_h(g)=\sum (Sh\o)gh\t$. Note also that
$\sum f^a\tens \Ad_{e_a}(h)$ is the left $\CA$-coadjoint coaction
on $\CH$ adjoint to the right adjoint coaction $\Ad(a)=\sum
a\t\tens (Sa\o)a\thr$  on $\CA$. \eproof

Next, since $\CA$ is the geometric quantity for us, it is useful
to recast these results in terms of an $\CA$-crossed module.
First, we can replace the left $\CH$ coaction on $W$ by a right
action of $\CA$ on $W$, and this by a left action of $\CA$ on
$W^*$. The result (after some computations) is \eqn{dcalca}{ \extd
h=\sum (h\t\la e_\beta\tens f^a\la f^\beta)\Ad_{e_a}(h\o)-\theta
h} \eqn{dcalcb}{ \extd a=\sum (e_\beta\tens \Ad_{f^a}(a\o)\la
f^\beta)e_a a\t-\theta a} \eqn{dcalcc}{ h.e_\alpha{}^\beta=\sum
(h\t\la e_\alpha\tens f^a\la f^\beta)\Ad_{e_a}(h\o)} \eqn{dcalcd}{
a.e_\alpha{}^\beta=\sum (e_\alpha\tens \Ad_{f^a}(a\o)\la
f^\beta)e_a a\t.} We can also write the left action of $\CH$ as a
right $\CA$-coaction, or if possible a left $\CA$-coaction on
$W^*$. Using the same notation for left coactions and $\sum
w\bz\tens w\bo$ for right coactions, the affected formulae become
\[ \extd h=\sum \<h\t, f^\beta\bo\>(e_\beta\tens f^a\la
f^\beta\bt)\Ad_{e_a}(h\o)\]
\[ h.e_\alpha{}^\beta=\sum \<h\t,e_\alpha\bo\>(e_\alpha\bz\tens f^a\la
f^\beta)\Ad_{e_a}(h\o).\]

Now suppose $\CA$ has a bicovariant calculus. Then $\bar\Lambda^1$
is a left $\CA$-crossed module as explained in Section~3 (a
quotient of $\CA^+$ under left multiplication and left adjoint
coaction). We therefore set $W=\bar\Lambda^{1*}$ in the above and
obtain an {\em induced calculus} on $D^*(\CA)$. One may check
easily that if the initial $\Omega^1(\CA)$ is inner with a
generator $\bar\theta$ say, then
\eqn{piOmegaD}{\pi(e_\alpha{}^\beta)=\<e_\alpha,\bar\theta\>
f^\beta,\quad \pi(h\tens a)=\eps(h)a} defines a surjection
$\Omega^1(D^*(\CA))\to \Omega^1(\CA)$ extending the canonical Hopf
algebra surjection $D^*(\CA)\to \CA$ as stated.

Of course not all crossed modules are of the type which comes
 from an initial differential structure on $\CA$. Let us give two examples,
one which is and one (the first) which is not. They both show a
different phenomenon whereby a calculus on $\CA$ (and on $\CH$)
can instead be sometimes included in the calculus on the double.
Note also that while the theory used is for finite-dimensional
Hopf algebras, all formulae can be used also with appropriate care
for infinite-dimensional Hopf algebras $\CA,\CH$ dually paired (or
skew-paired). For example, if $\CA=k[G]$ is an algebraic group of
Lie type dually paired with an enveloping algebra $\CH=U(\cg)$
then its double $D^*(G)=U(\cg)\lcocross k[G]$ is the tensor
product algebra and as a coalgebra crossed by the coadjoint
coaction, but the latter does not play a role in the differential
structure that results. We let $\{e_i\}$ be a basis of $\cg$. We
note that \eqn{vecfield}{ \del_if\equiv \sum\<e_i,f\o\>f\t,\quad
\forall f\in k[G]} is the action of the classical right-invariant
vector field generated by $e_i\in\cg$. We assume the summation
convention for repeated such indices. We use the 2-action
formulation (\ref{dcalca})-(\ref{dcalcd}).

\begin{example} We take $W=\bar\cg=1\oplus \cg\subset U(\cg)$ as a
subcrossed $U(\cg)$-module under the left adjoint action and the
regular coaction given by the coproduct. We write $e_0=1$ as the
additional basis element and let $\{f^i,f^0\}$ be the dual basis.
The left $U(\cg)$-crossed module structure and left action of
$k[G]$ on $W^*$ are
\[ \xi\la e_i=[\xi,e_i],\quad \xi\la e_0=0,\quad
\Delta_L e_i=e_i\tens e_0+1\tens e_i,\quad \Delta_Le_0=1\tens
e_0\]
\[ f\la f^i=f^i\, f(1),\quad f\la f^0=f^0 f(1)+ f^i\<e_i,f\>\]
for $\xi\in \cg$ and $f\in k[G]$. Then the resulting calculus
$\Omega^1(D^*(G))$ has structure
\[ \extd\xi= [\xi,e_i]\tens f^i+ e_0{}^i [\xi,e_i],\quad \extd f= e_0{}^i\del_i f \]
\[ [\xi,e_0{}^0]=e_0{}^i[\xi,e_i],\quad
[\xi,e_0{}^i]=0,\quad [\xi,e_i{}^0]=[\xi,e_i]\tens
f^0+e_i{}^j[\xi,e_j]\] \[ [\xi,e_i{}^j]=[\xi,e_i]\tens f^j,\quad
[f,e_\alpha{}^i]=0,\quad [f,e_\alpha{}^0]= e_\alpha{}^i\del_i f
;\quad \alpha=0,j.\]
\end{example}

Note that it is possible to restrict the calculus to basic forms
of the type $\{e^i{}_j,e_0{}^i\}$ as a subcalculus
$\Omega^1_{res}(D^*(G))$. When this is restricted to $k[G]$ we
have $\Omega^1_{res}|_{k[G]}=\Omega^1(G)$ its usual classical
calculus. When restricted to $U(\cg)$ we have
$\Omega^1_{res}|_{U(\cg)}$ the calculus with $\extd\xi=\rho(\xi)+e_0{}^i[\xi,e_i]$
where the first term is a standard type for $\Omega^1(U(\cg))$ as a
noncommutative space.  Thus $\Omega^1_{res}$
`factorizes' into something close to these standard constructions.
In fact we need not all $\span\{e_i{}^j\}$
here but only the image $\rho(U(\cg)^+)$ which could be different
unless $\cg$ is simple. Also note that while the classical
calculus on $k[G]$ is not inner, when viewed inside
$\Omega^1_{res}$ the element $\theta_0=\sum e_i{}^i$ generates
$\extd f$. The full calculus above is necessarily inner, by
construction.

As far as applications to physics are concerned let us note that
if one (perversely) regards $U(su_2)$ as a noncommutative $\R^3$,
i.e. as quantising the Kirillov-Kostant bracket then $D(SU_2)$ can
be viewed as the appropriate deformation of the isometry group of
$\R^3$. We refer to \cite{BatMa:non} for some recent work in this
area. Hence if one wants to construct an affine frame bundle etc.
as in Section~2 then one will need a calculus on $D^*(SU_2)$ such
as the above. Details will appear elsewhere \cite{BatMa:pre}.

Let us now give a different example on the same quantum double,
this time of the type induced by a calculus on $k[G]$. If one
tries to start with the classical calculus on $k[G]$, where
$\bar\Lambda^1=\cg^*$ the dual of the Lie algebra, one will get
zero for differentials $\extd f$ in the induced calculus. That is
why we had to work with an extension $\bar\cg$ above. Similarly,
at least in characteristic zero (and assuming an invertible
Killing form):

\begin{example} Let $\cg$ be a semisimple Lie algebra over $k$
and $c={1\over 2}K^{ij}e_i e_j$ the quadratic Casimir defined by
the inverse Killing form. We let $W=k c\oplus \cg=\tilde\cg\subset
U(\cg)^+$ as a subcrossed module under $\Ad$ and
$\Delta_L=\Delta-\id\tens 1$. We write $e_0=c$ to complete the
basis of $\cg$. The left $U(\cg)$-crossed module structure and
left action of $k[G]$ on $W^*$ are
\[ \xi\la e_i=[\xi,e_i],\quad \xi\la e_0=0,\quad
\Delta_L e_i=1\tens e_i,\quad \Delta_Le_0=1\tens
e_0+K^{ij}e_i\tens e_j\]
\[ f\la f^i=f^i\, f(1)+f^0 K^{mi}\<e_m,f\>,\quad f\la f^0=f^0 f(1)\]
for $\xi\in \cg$ and $f\in k[G]$. Then the resulting calculus
$\Omega^1(D^*(G))$ has structure
\[ \extd\xi=[\xi,e_i]\tens f^i+  e_i{}^0K^{mi}
[\xi,e_m],\quad \extd f= e_i{}^0K^{mi}\del_m f \]
\[ [\xi,e_0{}^0]=0,\quad [\xi,e_0{}^i]=e_0{}^0K^{mi}[\xi,e_m],\quad [\xi,e_i{}^0]=[\xi,e_i]\tens
f^0\] \[ [\xi,e_i{}^j]=[\xi,e_i]\tens f^j+e_i{}^0
K^{mj}[\xi,e_m],\quad [f,e_\alpha{}^0]=0,\quad
[f,e_\alpha{}^i]=e_\alpha{}^0 K^{mi}\del_m f.\]
\end{example}

We use the same conventions as the previous example. Again we have
a restricted subcalculus spanned by $\{e_i{}^j,e_i{}^0\}$ and
indeed the restriction of this to $k[G]$ is again the classical
calculus in terms of a new basis $e^m\equiv  K^{mi}e_i$. The
restriction to $U(\cg)$ is different, however. Unlike the previous
example, the input crossed module $\tilde\cg$ here is of the form
such that its dual is $\tilde\cg^*=\bar\Lambda^1$ for an initial
differential structure $\Omega^1_{init}(\CA)$ on $\CA$. Namely
\eqn{cascalc}{ \extd f=f^i\del_i f+f^0{1\over 2}K^{ij}\del_j\del_i
f,\quad [f,f^i]=f^0K^{mi}\del_m f,\quad [f,f^0]=0} for all $f\in
k[G]$. This is a standard 1-dimensional noncommutative extension
of the usual classical calculus in which the second-order
Laplacian evident here is viewed as a `first order' partial
derivative $\del_0$ in the extra $f^0$ direction (the other
$\del_i$ are the usual classical differentials as above). Thus a
reasonable but noncommutative calculus can induce a reasonable one
on the double.

Finally, we consider the above results as a step towards a
`T-duality' theory for differential calculi whereby Hopf algebra
duality is extended to noncommutative geometry. In physics,
Poisson-Lie T-duality refers to an equivalence between a
$\sigma$-model on Poisson Lie group $G$ and one on its Drinfeld
dual $G^*$ with dual Lie bialgebra, see
\cite{KliSev:dua,BegMa:poi} and elsewhere. The transfer of
solutions is via the Lie bialgebra double $D(\cg)$. We could hope
for a similar theory in the quantum group case. For this one would
need first of all to extend the Hopf algebra duality functor to
differentials. Our construction indicates a way to do this: given
a calculus on $\CA$ we can induce one on the codouble and then
project that down to one on $\CH$. Probably some modifications of
this idea will be needed for nontrivial results (as the examples
above already indicated) but this is a general idea for
`transference of calculi' that we propose.

\bigskip
\subsection*{Acknowledgements} I would like to thank Xavier Gomez
for discussions. The author is a Royal Society University Research
Fellow.


\end{document}